\documentstyle[12pt]{article}
\def\ben{\begin{enumerate}}
\def\een{\end{enumerate}}
\def\bel{\begin{equation}\label}
\def\ee{\end{equation}}
\def\ba{\begin{array}}
\def\ea{\end{array}}
\def\a{\alpha}
\def\r#1{(\ref{#1})}
\def\cite#1{[#1]}
\def\i{\item}
\def\tens{\otimes}
\newcounter{popnr}
\def\fn#1{{\mathop{{\rm #1}}}}
\def\theequation{\thesection.\arabic{equation}}
\renewcommand{\theequation}{\arabic{section}.\arabic{equation}}
\newcommand{\alpheqn}{\setcounter{popnr}{\value{equation}}
                      \addtocounter{popnr}{1}
                      \setcounter{equation}{0}
   \renewcommand{\theequation}{\arabic{section}.
  \arabic{popnr}\alph{equation}}}
\newcommand{\reseteqn}{\setcounter{equation}{\value{popnr}}
     \renewcommand{\theequation}
     {\arabic{section}.\arabic{equation}}}
\def\bl{\alpheqn}
\def\el{\reseteqn}
\def\sec{\setcounter{equation}{0}}

\def\~{\widetilde}
\def\cop{\Delta}
\def\ftm{\footnotemark}
\def\to{\,\longrightarrow\,}
\def\lam{\lambda}
\def\k{\kappa}
\def\poin{Poincar\'e }
\begin{document}
\begin{titlepage}
\title{Quantum Deformations of Conformal Algebras Introducing
Fundamental
Mass Parameters}
\author{Jerzy Lukierski \thanks{On leave of absence from the
Institute of
Theoretical Physics, University of Wroc{\l}aw, 50-204 Wroc{\l}aw,
Poland.}
Scientific
\thanks{Partially supported by KBN grant 2P 302 087 06.},
Pierre Minnaert
and Marek Mozrzymas \ftm[1]\\
Laboratoire de Physique Th{\'e}orique,\\ Universit{\'e} Bordeaux I,\\
 33175 Gradignan, France}
\date{}
\maketitle
\begin{abstract}
We consider new class of classical $r$-matrices for $D=3$ and $D=4$
conformal Lie  algebras.
These r-matrices do satisfy the classical Yang-Baxter equation and as
two-tensors
belong to the tensor
product of Borel subalgebra. In such a way we generalize
the lowest order of known
nonstandard quantum
deformation of $sl(2)$ to the Lie algebras $sp(4) \cong so(5)$ and
$sl(4) \cong so(6)$.
As an exercise we interpret nonstandard deformation of $sl(2)$ as
describing quantum $D=1$ conformal algebra with fundamental mass
parameter. Further we describe the $D=3$ and $D=4$ conformal
bialgebras
with deformation parameters  equal to the inverse of fundamental
masses. It appears that for $D=4$ the deformation of the \poin algebra
sector coincides with "null plane" quantum \poin algebra.
\end{abstract}
\vfill
CPTMB/PT/95-6
\hfill
June 1995
\end{titlepage}

\section{Introduction}
Recently, quantum deformations of $D=4$ space-time symmetries were
considered (see e.g. [1--16]). Because the \poin generators ($P_\mu$,
$M_{\mu\nu}$) transform nontrivially under the change of length
scale
\bel{1.1}
P_{m}^{\prime}=\lambda^{-1}P_{\mu}\,,\qquad M_{\mu
\nu}^{\prime}=M_{\mu \nu}\,,
\ee
in the case of $D=4$ \poin
algebra one can distinguish two different types of quantum
deformations:
\ben
\i[a)]
With dimensionless deformation parameter $q$ \cite{2,4,12,15}.

In such a case the deformed \poin algebra $U_q(P_4)$ is invariant
under the
rescaling \r{1.1} with  the value of $q$ kept fixed.
\i[b)]
With dimensionfull mass-like parameter $\k$ \cite{1,3-5,7-9,13,14}.

The $\kappa$-deformed \poin algebra, denoted by $U_{\kappa}(P_4)$ is
invariant
under rescaling (1.1) provided that we rescale also $\kappa
\longrightarrow
\kappa^{\prime}=\lambda^{-1} \kappa$.
Such a deformation occurs naturally in the  bicrossproduct
description of
quantum
\poin algebra \cite{8}, with the deformed semidirect product of the
Lorentz
algebra and fourmomentum sector containing nonlinear functions of the
fourmomenta. The  scale invariance of these nonlinear functions
under the transformation (1.1) leads
necessarily to the appearence of a mass-like deformation parameter.
\een

In this paper we shall consider the quantum deformations of the
conformal
algebras  with a dimensionfull deformation parameter.
The importance of this problem can be justified as follows:
\ben
\i[i)]
The conformal symmetry is the fundamental ``master'' symmetry of
space-time,
containing two other fundamental geometries (\poin, de Sitter) as its
broken
cases.
\i[ii)]
The conformal symmetries describe the world of massless particles and
fields. Usually, masses breaking conformal symmetries are introduced
on
the level of field representations of the symmetry group. The
introduction
of the mass-like deformation parameter leads to the appearence of the
fundamental mass on basic
geometrical level.
\i[iii)]
{}From the mathematical point of
view, the deformations of conformal algebras with fundamental mass
parameter
introduce a new type of quantum deformations of Lie algebras,
generalizing
the nonstandard deformation of $sl(2)$ \cite{17-21}.
\een

The quantum deformations of Lie algebras are described by the
noncommutative
and noncocommutative Hopf algebras, which can be obtained by the
quantization
of
corresponding Lie bialgebras \cite{22, 23}. For semisimple Lie
algebras the
classification of quantum deformations is provided by the choice of
classical $r$-matrices, determining the cocommutators of Lie
bialgebra. In
particular the classical Yang-Baxter equations is required if we
assume
that after
completing the quantum deformation procedure we get the
quasitriangular
Hopf algebra with
quantum universal $R$-matrix, satisfying quantum Yang-Baxter equation
\cite{17}.

In this paper we shall consider the choices of classical $r$-matrices,
generalizing for  $sp(4) \cong so(5)$ and $sl(4) \cong so(6)$ the
nonstandard
deformation of $sl(2)$
\cite{17-21},
described by the classical $r$-matrices $r_\pm =h\wedge e_\pm$
($A\wedge
B\equiv (A\tens B - B \tens A)$).

The plan of this paper is the following:

In Sect.\ 2 we recall a known example of nonstandard deformation of
$sl(2)$ Lie algebra and we interprete it as the quantum deformation
od $D=1$
conformal algebra introducing mass-like deformation parameter.

In Sect.\ 3 we  introduce for any simple Lie algebra $\hat{g}$ the
class of nonstandard classical $r$-matrices described by two-tensors
on its
Borel subalgebra $\hat{b}$. In particular we shall be interested in
maximal
nonstandard classical $r$-matrices of $\hat{g}$, i.e. those which
cannot be
described by the generators of any simple subalgebra
$\hat{g}^{\prime} \subset
\hat{g}$ ($\hat{g}^{\prime} \neq \hat{g}$). We present  these
maximal nonstandard classical $r$-matrices in Cartan-Weyl basis for
the complex
Lie algebras $so(5) \cong sp(4)$ and $so(6) \cong sl(4)$; further we
shall
impose the $so(3,2)$ and $so(4,2)$ reality conditions.

In Sect.\ 4 and \ 5 we introduce the physical (conformal) generators
for $so(3,2)$
($D=3$ conformal) and $so(4,2)$ ($D=4$ conformal). It appears that the
nonstandard $r$-matrices considered in sec.\ 3 lead to the
introduction
of the
deformation parameters $\frac{1}{M_i}$ ($M_i$ - fundamental masses).
More specifically we obtain that
\ben
\i[i)]
For $D=3$ the deformation is described by two masses $M_1$, $M_2$. If
$M_2=0$ our classical $r$-matrix generates the quantum defomation of
$D=3$
\poin algebra.
\i[ii)]
For $D=4$ conformal algebra the deformation generated by our maximal
nonstandard $r$-matrix introduces one fundamental mass $M$ and
describes the
defomation
of the \poin subalgebra with classical fourdimensional subalgebra
$e(2) \oplus r$ ($e(2)$ is $D=2$ inhomogenous Euclidean algebra and
$r$ generator of noncompact Abelian group $R$).
It is interesting that such a deformation has been obtained recently
in
\cite{24} under the name of "null-plane" quantum  \poin algebra
in different context - by so-called
deformation embedding method \cite{25} from three to four dimensions.
\een
In Sect.\ 6 we present the discussion and outlook.

\sec
\section{Quantum deformations of $D=1$ conformal algebra with mass-
like
deformation parameter}

It is well known that the $D=1$ conformal algebra
\bel{2.1}
[D,P]=P\,,\qquad [D,K]=-K\,,\qquad[P,K]=2D\,,
\ee
where $P$ describes the time translations (energy), $D$ --- the
scaling
transformations, and $K$ --- the conformal accelerations, after the
identification
\bel{2.2}
e_+=P\,,\qquad e_-=K\,,\qquad h=2D\,,
\ee
can be written as the $sl(2,R)\simeq O(2,1)$ algebra in Cartan-
Chevalley
basis
\bel{2.3}
[h,e_\pm]=\pm 2 e_\pm\,,\qquad [e_+,e_-]=h\,.
\ee
For $sl(2,R)$ there exist only two inequivalent deformations,
generated by
the classical $r$-matrices with  the following antisymmetric parts:\\
i)standard deformation \cite{22--24}
\bel{2.4}
r_s=c_s e_+ \wedge e_-\,.
\ee
Adding suitable symmetric part one gets from (2.4) the linear term in
the
coproduct of
Drinfeld-Jimbo quantum algebra $U_q(sl(2))$. It is easy to see that
the
invariance of (2.4) under the rescalings $ P \to \lam^{-1}P$,
$K\to\lam K $
imply that the deformation parameter $c_s$ is dimensionless.\\
ii)nonstandard deformations \cite{17-21}.
\bel{2.5}
r_\pm=c_\pm h \wedge e_\pm\qquad (c_\pm\mbox{ --- constants}).
\ee
Using (2.2) one obtains
\bel{2.6}
r_+ =2 c_+ D \wedge P\,,\qquad r_-= 2c_- D\wedge K\,.
\ee
We see that the quantum deformation generated by the classical $r$-
matrix
$r_+$ provides the deformation parameter transforming as the inverse
of mass
($2c_+=\frac1M$ ; we shall call it $M$-deformation) and $r_-$ implies
that
$2c_- =\~M$ ($\~M$-deformation), where $M$, $\~M$ are fundamental
masses.\\
The quantization of the Lie algebra $sl(2,R)$ generated by $r_+$ has
been
given firstly by Ohn \cite{21}. The relations (2.1) in all orders of
$\frac{1}{M}$ are deformed as follows:
\bel{2.7}
\ba{rcl}
[D,P]&=&M \sinh \frac P M\,, \qquad [P,K]\,=\, 2D\,,\\[2mm]{}
[D,K]&=& \frac{1}{2}(-K \cosh \frac P M - (\cosh \frac P M )K)\,.
\ea
\ee
The coproduct and the antipode take the form:
\bel{2.8}
\ba{rcl}
\cop (P) &=& P \tens 1 + 1 \tens P\,, \\
\cop (D) &=& e^{-\frac P M} \tens D + D \tens e ^{\frac P M}\,,\\
\cop (K) &=& e^{-\frac P M} \tens K + K \tens e ^{\frac P M}\,,
\ea
\ee
\bel{2.9}
S(P)=-P ,\qquad S(K)\,=-K-\frac{1}{M}(D-\sinh\frac{P}{M})\,
\qquad S(D)\,=-D-2\sinh\frac{P}{M}.
\ee
We would like to make  the following remarks:
\ben
\i[i)]
The quantum deformation (2.7)-(2.8) has the $D=1$ quantum Weyl Hopf
subalgebra span by two generators $P$ and $D$.
\i[ii)]
The deformation generated by the classical $r$-matrix $r_-$ (see
(2.5)-(2.6)) is obtained by the following replacement in the formulae
(2.7)-(2.9)
\bel{2.10}
\ba{c}
P\to K\,,\qquad K\to P \,,\qquad D\to -D \,,\\
M\to \frac 1 {\~M} \,.
\ea
\ee
It appears that the $M$-deformation is described by the  power series
in
the variable $M^{-1}$ (like the $\k$-deformation of \poin algebra)
and the
$\~M$-deformation --- as the  power series in $\~M$. Respectively, the
``classical''
no-deformation limits are $M\to \infty$ and $\~M\to 0$.
\i[iii)]
Because  $D=2$ conformal algebra $so(2,2)\simeq so(2,1)\oplus
so(2,1)$, the
deformation of $D=2$ conformal algebra is determined by the pair of
(standard or
nonstandard) deformations of $so(2,1)\cong sl(2,R)$. In fact the
decomposition of the
conformal algebra into one-dimensional conformal algebras of right-
movers
and left-movers permits to introduce different $M$  and $\~M$
parameters in
these two chiral sectors.
The description of quantum $D=2$ conformal algebra by a pair of
nonstandard
deformations with one parameter was proposed recently in \cite{26}.
The application of
nonstandard deformation of $so(2,1)$ to the description of
quantum-mechanical model is also under consideration \cite{27}.\\
\een

\sec
\section{Nonstandard classical $r$-matrices for $sp(4;C)$, $sl(4;C)$
and
their real forms}

We shall consider the
nonstandard deformations of a simple Lie algebra $\hat g$
generated by classical $r_\pm$ -matrices satisfying the following
two conditions:
\ben
\i[a)] $r_\pm \in \hat b_\pm \wedge \hat b_\pm$, where  $\hat b_\pm$
describes the Borel subalgebras, span respectively by ($h_i$, $e_{\pm
a}$)
($i=1\ldots r=\fn{rank } \hat g$, $a=1\,\ldots N=\frac12 (\fn{dim}
\hat g
-r)$), where $h_i$ describes the Abelian Cartan generators, and
$e_{+a}
(e_{-a})$the positive (negative) root generators.
\i[b)]
$r_\pm$ satisfy the classical Yang-Baxter equation
\bel{3.1n}
[r_{\pm12}, r_{\pm13}] + [r_{\pm12},r_{\pm23}]+
[r_{\pm13},r_{\pm23}]=0\,,
\ee
where if $r=r^{AB} I_A\otimes I_B$  ($I_A$ -- Cartan-Weyl basis of
$\hat g$;
$A,B=1\ldots \fn{dim} g$), $r_{12}= r^{AB}I_A\tens I_B \tens 1$ etc.
\een
We shall introduce also the notion of maximal nonstandard $r$-matrix,
belonging
to the class described above. Let us observe that if
$\hat{g}^{\prime}\subset
\hat{g}$ ($\hat{g}^{\prime}$ is a simple Lie subalgebra of
$\hat{g}$), then
any classical $r$-matrix for $\hat{g}^{\prime}$ is also a
classical $r$-matrix for $\hat{g}$. For the maximal nonstandard
classical
$r$-matrices
such an embedding does not exist - the algebra $\hat{g}$ is the
minimal simple
algebra,
providing Borel algebras $\hat{b}_{\pm}$ in the formula $r_{\pm}\in
\hat{b}_{\pm}
\wedge\hat{b}_{\pm}$.
Below we shall consider only the  maximal nonstandard classical $r$-
matrices
$r_{+}\in \hat{b}_{+} \wedge\hat{b}_{+}$. It appears that the matrices
$r_{-}\in \hat{b}_{-} \wedge\hat{b}_{-}$ can be obtained by suitable
anti-
automorphism of $\hat{g}$.\\
In physical applications there are relevant real Lie algebras.
We shall consider therefore only
the  classical $r$-matrices satisfying the reality condition:
\bel{3.2}
(r_+)^{+} = ((r_+)^{AB}I_A \wedge I_B)^{+} =
((r_+)^{AB})^{\ast}I_A^{+}
 \wedge I_B^{+} = r_+
\ee
We assume that $+$--involution is an anti-automorphism of $U(\hat{g})$
(i.e. $(ab)^{+} = b^{+}a^{+}$) and automorphism of tensor product
(i.e. $(a\otimes b)^{+} = a^{+}\otimes b^{+}$). Because from (3.2)
follows that $(b_+)^{+} \subset
b_+$, the involutions implying (3.2) map $\Delta^{\pm} \longrightarrow
\Delta^{\pm}$ where by $\Delta^+ (\Delta^-)$ we denote the set of
positive
(negative ) root generators. It appears that due to the parabolic
decomposition of the conformal algebra $so(D,2)$ (see e.g. \cite{4})
indeed such a type of reality conditions describe real conformal
algebras.\\
The general framework describing solutions of the classical Yang-
Baxter
equation, related to the Borel subalgebras of simple Lie algebras was
considered in \cite{28} (see also \cite{29}).
Below we shall describe two examples:
\ben
\i[a)]$Sp(4,C)\cong so(5,C)$ ($D=3$ complex conformal algebra)\\
The Borel subalgebra $b_+$ has the basis ($h_1, h_2, e_1, e_2, e_3,
e_4$)
where $e_1, e_2, $ are the simple root generators (Cartan-Chevalley
basis)
and
\bel{3.3}
e_3=[e_1, e_2],\qquad e_4=[e_1, e_3],
\ee
The maximal nonstandard classical $r$-matrix takes the form
\bel{3.4}
r_+ = c_3^{(1)}(h_1\wedge e_4 - e_1 \wedge e_3)+ c_3^{(2)}h_2\wedge
e_4
\ee
The $+$-involution introducing the real form $Sp(4,R)\cong so(3,2)$ we
choose the following (see e. g. \cite{30}; $i=1,2$):
\bel{3.5}
\ba{r@{=}lr@{=}l}
h_i^{+}&-h_i\,,\\[2mm]
e_1^{+}&\lambda e_1\,,&\qquad e_2^{+}&\epsilon e_2\,,\\[2mm]
e_3^{+}&-\lambda \epsilon e_3\,,&\qquad e_4^{+}&\epsilon e_4\,,
\ea
\ee
where $\lambda^2 = \epsilon^2 =1$.\\
{}From the invariance of (3.4) under (3.5) follows that
$\epsilon=\pm1, \lambda = \pm1$. One gets for real $c_3^{(1)}$,
$c_3^{(2)}$ that $\epsilon=-1$, $\lambda =\pm 1$.
It can be seen from \cite{31} that both real forms describe the real
conformal
algebra $so(3,2)$.
\i[b)] $sl(4)\cong so(6,C)$ ($D=4$ complex conformal algebra).\\
The simple root generators $e_{\pm i}$
($i=1,2,3$) and Cartan generators $h_i$ define the remaining part of
Cartan-Weyl basis as follows:
\bel{3.6}
\ba{r@{=}lr@{=}l}
e_4&[e_1,e_2]\,,&\qquad e_{-4}&[e_{-2},e_{-1}]\,,\\[2mm]
e_5&[e_2,e_3]\,,&\qquad e_{-5}&[e_{-3},e_{-2}]\,,\\[2mm]
e_6&[e_1,e_5]\,,&\qquad e_{-6}&[e_{-5},e_{-1}]\,,
\ea
\ee
and
\bel{3.7}
h_4=h_1+h_2\,,\qquad h_5=h_2+h_3\,, \qquad h_6=h_1+h_2+h_3\,.
\ee
The $sl(4)$ Lie algebra can be written compactly as follows
($A,B=1,\,\ldots 6$)
\bel{3.8}
\ba{rcl}
[h_A,e_{\pm B} ] &=& \pm \a_{AB} e_{\pm B} \qquad\mbox{(no
summation)}\,,\\[2mm]{}
[e_A,e_{-B} ] &=& \delta_{AB} h_B \qquad\mbox{(no summation)}\,,
\ea
\ee
where the extended symmetric Cartan matrix is given by the formula
(see e.g.\ \cite{32})
\bel{3.9}
\a_{AB}=\pmatrix{       2 &-1&0 &1 &-1&1\cr
                        -1&2 &-1&1 &1 &0\cr
                        0 &-1&2 &-1&1 &1\cr
                        1 &1 &-1&2 &0 &1\cr
                        -1&1 &1 &0 &2 &1\cr
                        1 &0 &1 &1 &1 &2}
\ee
The remaining $sl(4)$ Lie algebra relations (besides the ones given by
(3.6)) are obtained from Serre relations for the generators $e_{\pm
i}$
($i=1,2,3$) and take the form:
\bel{3.10}
\ba{lllllllll}
[e_1,e_3]&=&[e_1,e_4]&=& [e_1,e_6]&=&0\,,\\[2mm]{}
[e_2,e_4]&=&[e_2,e_5]&=& [e_2,e_6]&=&0\,,\\[2mm]{}
[e_3,e_4]&=&e_6\,,   & &  [e_3,e_5]&=&[e_3,e_6]&=&0\,,\\[2mm]{}
[e_4,e_5]&=&[e_4,e_6]&=& [e_5,e_6]&=&0\,.
\ea
\ee
\een
For $sl(4,C)$ we found the following maximal nonstandard classical r-
matrix
\bel{3.11}
r_+ = c_4^{(1)}(h_1 - h_3)\wedge e_6 + c_4^{(2)}(h_3\wedge e_6 +
e_1 \wedge e_5 - e_3 \wedge e_4)
\ee
{}From the discussion of all real forms for $sl(4,C)$ (see e.g.
\cite{32,33}
and put $q=1$) we select
one which maps $\Delta^{\pm} \longrightarrow \Delta^{\pm}$  and
describes the
conformal algebra $so(4,2)$:
\bel{3.12}
\ba{r@{=}lr@{=}l}
h_1^{+}&-h_3\,,&\qquad h_2^{+}&-h_2\,,\\[2mm]
e_1^{+}&\epsilon e_3\,,&\qquad e_2^{+}&\eta e_2\,,\\[2mm]
e_4^{+}&\eta \epsilon e_5\,,&\qquad e_6^{+}&\eta e_6\,,
\ea
\ee
{}From (3.11-12) follows that the reality conditions (3.12) imply
$\eta =\pm 1, \epsilon = \pm1$. For $\eta =-1$, $\epsilon =\pm 1$
we choose  $c_4=c^{(2)}_4=2c^{(1)}_4$ ($c_4$ real).

\sec
\section{$D=3$ Conformal Algebra }

For the description of $D=3$ conformal algebra in terms of
$Sp(4)\cong so(5)$ Cartan-Weyl basis satisfying the reality
conditions (3.5)
we shall use the formulae given in \cite{30} (we put for them
$q=1$). Selecting for simplicity $\lambda = -\epsilon =1$ we obtain
\bl
\bel{4.1a}
\ba{rclrcl}
h_1 &=& M_{12}\,,\qquad& h_2&=& M_{04}-M_{12}\,,\\[2mm]
e_1 &=& \frac{1}{\sqrt{2}}(M_{23}+M_{31})\,,\qquad &e_3 &=&
\frac{1}{\sqrt{2}}(M_{03}+M_{34})\,
\ea
\ee
\bel{4.1b}
\ba{rclrcl}
e_2 &=& -\frac{1}{\sqrt{2}}(M_{14}+M_{24}+M_{01}+M_{02})\,,\\[2mm]
e_4 &=& \frac{1}{\sqrt{2}}(M_{14}+M_{01}-M_{24}-M_{02})\,.
\ea
\ee
\el
where $M_{AB}^{\dagger}=-M_{AB}=M_{BA}$ ($A,B=0, 1, 2, 3, 4 $), and
it follows
that
\cite{30,31}
\bel{4.2}
[M_{AB},M_{CD}] = \eta_{BC}M_{AD} +\eta_{AD}M_{BC}-\eta_{AC}M_{BD}
-\eta_{BD}M_{AC}\,,
\ee
where $\eta_{AB}=\fn{diag} (-1,1,-1,1,1,)$. The generators $M_{31}= J,
M_{12}=L_1, M_{23}=L_2$ form $D=3$ Lorentz algebra $so(2,1)$. The
rest of
conformal generators are defined as follows:
\bl
\bel{4.3.a}
\ba{r@{=}lr@{=}l}
M_{01}& \frac{1}{\sqrt{2}}(P_1+K_1)\,,
&\qquad M_{41}& -\frac{1}{\sqrt{2}}(P_1-K_1)\,,\\[2mm]
M_{02}& \frac{1}{\sqrt{2}}(P_0+K_0)\,,
&\qquad M_{42}& -\frac{1}{\sqrt{2}}(P_0-K_0)\,,\\[2mm]
M_{03}& \frac{1}{\sqrt{2}}(P_2+K_0)\,,
&\qquad M_{43}& -\frac{1}{\sqrt{2}}(P_2-K_0)\,,
\ea
\ee
\bel{3.12c}
D = M_{04}\,.
\ee
\el
Substituting this physical   basis into the formulae (4.1) one obtains
that ($P_{\pm} = P_0 \pm P_1$)
\bel{4.4}
r_+ =\frac{1}{M_1}(L_1\wedge P_1 - (L_2+J) \wedge P_2)+
\frac{1}{M_2}(D-L_1)\wedge P_-
\ee
where$M_1=2(c^{(1)}_3)^{-1}, M_2=2(c^{(2)}_3)^{-1} $ describe the
fundamental mass
parameters. One can make the following remarks:
\ben
\i[i)]If $M_2= \infty, M_1< \infty$, the classical $r$-matrix (4.4)
describes the deformation of  $D=3$ \poin algebra. This deformation
is different from the $\kappa$-deformation of $D=3$ \poin algebra
described by the following  classical  $r$-matrix \cite{13}
\bel{4.5}
r_+ =\frac{1}{\kappa}(L_1\wedge P_1 + L_2 \wedge P_2)
\ee
\i[ii)]Two other selected cases are $M_1=M_2$ and $M_1=\infty$. The
second case describes due to relation $[D-L_1, P_-]=0$  so called
soft deformations \cite{34}.
\een

\sec
\section{$D=4$ Conformal Algebra}
We shall express the $sl(4,C)$ Cartan-Weyl generators satisfying the
reality
conditions (3.12) by using the formulae given in \cite{32}. The
relations
(3.8) and (3.10) imply $so(4,2)$ classical $D=4$ conformal algebra
relations ($M_{KL}^{\dagger}=-M_{KL}=M_{LK}$; $K,L=$$0, 1, 2, 3, 4,
5$)
\bel{5.1}
[M_{KL},M_{PR}] = g_{KR}M_{LP} +g_{LP}M_{KR}-g_{LR}M_{KP}
-g_{KP}M_{LR}\,,
\ee
where $g_{AB}=\fn{diag} (-1,1,1,1,1,-1)$ and
\bel{3.11}
\ba{rclrclrcl}
P_\mu&=&(M_{4\mu}+M_{5\mu})\,,\qquad&K_\mu&=&
(M_{5\mu}-M_{4\mu})\,,\\[2mm]
M_i&=& \frac12 \epsilon_{ijk}M_{jk}\,,& L_i
&=&M_{0i}\,,\qquad&D&=&M_{45}\,
\ea
\ee
Introducing $M_{\pm}=M_1\pm iM_2, L_{\pm}=L_1\pm iL_2$ one gets
(we choose $\epsilon=1$ in (3.12))
\bl
\bel{5.3a}
\ba{r@{=}lr@{=}l}
h_1& L_3 -iM_3)\,, &\qquad h_3 & L_3+iM_3 \,,\\[2mm]
e_1& \frac12(M_+ +iL_+)\,, & e_3 &-\frac12(M_- -iL_-)\,,\\[2mm]
e_2& \frac{1}{2}(P_0 -P_3)\,, & e_6 &\frac{1}{2}(P_0+ P_3)\,,\\[2mm]
e_4& \frac{i}{2}(P_1 +iP_2)\,, & e_5 &-\frac{i}{2}(P_1- iP_2)\,,
\ea
\ee
\bel{5.3b}
h_2 = -(D+L_3)\,.
\ee
\el
Substituting this physical basis into the formula (3.11) with
$c_4=c_4^{(2)}=2c_4^{(1)} =\frac{2}{M}$ one obtains ($P_+=P_0+P_3$)
\bel{5.4}
r_+ = \frac{1}{M}[L_3\wedge P_+ +(M_2 + L_1)\wedge P_1
-(M_1 - L_2)\wedge P_2]
\ee
The algebra $\~{e}(2)$ with the generators
\bel{5.5}
\~{E}_1 = L_1 +M_2,\qquad \~{E}_2 = -L_2 +M_1,\qquad \~{E}_3 = L_3,
\ee
has the following commutation relations
\bel{5.6}
[\~{E}_1,\~{E}_2] = 0,\qquad [\~{E}_1,\~{E}_3]=-\~{E}_1,\qquad
[\~{E}_2,\~{E}_3]=-\~{E}_2
\ee
and describes $D=2$ \poin algebra.\\
One can decompose the Lorentz algebra $so(3,1)$ as follows
\bel{5.7}
so(3,1)=\~{e}(2)\oplus e(2)
\ee
where $e(2)$ generators are the following
\bel{5.8}
E_1 = L_1 -M_2,\qquad E_2 = L_2 +M_1,\qquad E_3 = M_3,
\ee
and satisfy the $D=2$ Euclidean algebra relations
\bel{5.9}
[E_1,E_2] = 0,\qquad [E_2,E_3]=-E_1,\qquad
[E_1,E_3]=E_2,
\ee
The M-deformation, generated by the classical $r$-matrix (5.4) leads
to the
quantum deformation of the \poin algebra obtained recently in
\cite{24} by
different method. It is easy to check that the classical $r$-matrix
(5.4)
does not modify the coproducts for the generators ($M_3, \~{E}_1,
\~{E}_2$),
forming another $D=2$ Euclidean algebra as well as
the component $P_+$ ($P_+ = P_0 +P_3$) of the four-momentum.
The fourdimensional algebra with the generators ($M_3, \~{E}_1,
\~{E}_2,
P_+$) describe the classical subalgebra of the $M$-deformed $D=4$
coformal algebra.
The mass
Casimir is given by the formula \cite{24}
\bel{5.10}
C_2 = P_1^2 + P_2^2 -M P_-sinh(\frac{P_+}{M})
\ee
Because the deformation $U_M({\cal P}_4)$ of the \poin algebra, called
in \cite{24} the "null plane" quantum \poin algebra forms a Hopf
subalgebra of $U_M(o(4,2))$,
it is the deformed mass-shell (5.10) which describes the
deformation of massless d'Alembert operator in a way which permits
also
to introduce the $\kappa$-deformed conformal properties. Having the
formula for
the classical $r$-matrix one can look for the quantum deformations of
the
generators $K_{\mu}$ and $D$ using e.g. the perturbative formulae for
quantisation of bialgebras  given by Drinfeld \cite{35} and
Reshetikhin
\cite{36}. We conjecture that similarly as the relation $\vec{P}^2-
P_0^2=0$
is conformal invariant in undeformed theory, we shall obtain that
for $M<\infty$
\bel{5.11}
[C_2,K_{\mu}]\mid_{C_2=0}=[C_2,D]\mid_{C_2=0} =0
\ee
The proof of  relation (5.11) would additionally justify
the choice of the deformation  generated by classical $r$-matrix (5.4)
as the most
appropriate for the massless conformal-invariant theories and explain
the
difficulties with embedding of $\kappa$-\poin algebra
$U_{\kappa}({\cal P}_4)$ into the quantum conformal algebra
(see \cite{37}).\\
\sec
\section{Discussion and Outlook}
In this paper we proposed new deformation schemes for $D=3$ and $D=4$
conformal algebras introducing fundamental mass parameters. The
existence of a fundamental scale as a lower bound to any position
measurement seems to be a consequence of quantisation of general
relativity (see e. g. \cite{38}). The considerations of this paper
introduce on purely geometric basis the notion of deformed conformal
structures, which should modify the classical notions of distance
and causality for short distances. By suitable adjustment of
the mass-like  deformation parameters one should be able to restrict
such a modifications to the distances comparable or shorter than
the Planck length ($\Delta l \simeq 10^{-33}cm$).\\
The deformations described by classical $r$-matrices (4.4) and (5.4)
describe new $D=3$ and $D=4$ conformal bialgebras. If we introduce
the fundamental matrix realizations of $D=3$ and $D=4$ classical
conformal algebras ($4\times 4$ real matrices for $D=3$ and
$4\times 4$ complex matrices for $D=4$) one can write down
the Poisson brackets describing the Lie-Poisson structure on
$D=3$ and $D=4$ conformal groups \cite{22}.Following the
derivation in \cite{13} of $D=4$ quantum \poin group one can
quantize the Poisson-Lie brackets and obtain quantum $D=3$
and $D=4$ conformal groups.\\

We would like to recall here that the $\kappa$-deformation of \poin
algebra (see \cite{1, 3, 8}) leads to classical nonrelativistic
symmetries $so(3)\subset so(3,1)$  but the "null plane" quantum
\poin algebra given in \cite{23} provides classical symmetries
for the subgroup $e(2)\subset so(3,1)$. The subgroups $so(3)$ and
$e(2)$ are respectively the stability groups for the time-like and
light-like four-momenta. One can also introduce the third deformation,
leading to the classical symmetries for the subalgebra $so(2,1)\subset
so(3,1)$ not changing the components of the space-like (tachyonic)
four-momenta, by assuming the following classical $r$-matrix
(see \cite{39})
\bel{6.1}
r=M_1\wedge P_2-M_2\wedge P_1+L_3\wedge P_0
\ee
where we have chosen as the "deformation direction" the third space
axis. It is quite possible that \underline{different} quantum
deformations of the $D=4$ \poin algebra are appropriate for the
world with matter propagating with sublumimal velocities,
light velocities and  superluminal  velocites. These three
cases can be considered in a nice geometric setting if we
introduce the quantum deformation of the \poin algebra as  the
bicrossproduct Hopf algebra \cite{8}:
\bel{6.2}
U({\cal P}_4) = U(so(3,1)\bowtie T_4^{\kappa}
\ee
The three classes of deformations discussed above are defined
by the choice in (6.2) of the
part of the Lorentz algebra $so(3,1)$ remaining unaffected by
the $\kappa$-deformed action of $so(3,1)$ generators on
the four-momentum sector $T_4^{\kappa}$.
This way of looking at different deformations of \poin algebra is
under
considerations.\\

{\bf Acknowledgement}\\
One of the author (J. L.) would like to thank C. Gomez, O. Ogievetsky,
H. Ruegg,
M. Santander, A. Stolin and V. Tolstoy for discussions. Two authors
(J.L. and M.M.)
would like to
thank the University of Bordeaux I for the hospitality and financial
support.\\

\frenchspacing
\def\pl#1{{\it  Phys. Lett. } {\bf #1}}
\def\comm#1{{\it  Phys. Lett. } {\bf #1}}
\def\jp#1{{\it Journ. Phys.} {\bf #1}}
\def\jmp#1{{\it Journ. Math. Phys.} {\bf #1}}
\def\lmp#1{{\it Lett. Math. Phys.} {\bf #1}}
\def\bib{\bibitem}

\end{document}